\documentclass{article}
\usepackage[ansinew]{inputenc}
\usepackage[dvips]{graphicx}
\usepackage[T1]{fontenc}
\usepackage{amsmath,enumerate,amssymb,hhline,amsthm}

 \newtheorem{thm}{Theorem}[section]
 \newtheorem{cor}[thm]{Corollary}
 \newtheorem{lem}[thm]{Lemma}
 
 \newtheorem{conj}[thm]{Conjecture} 
 \theoremstyle{definition}
 \newtheorem{rem}[thm]{Remark}

\begin{document}

\title{On the binary additive divisor problem in mean}

\author{Eeva Suvitie}

\date{}

\maketitle

\let\thefootnote\relax\footnotetext{This
research was supported by the Finnish Cultural Foundation and by the project 138522 of the Academy of Finland, and
accomplished while the author was visiting \'{E}cole Polytechnique F\'{e}d\'{e}rale de Lausanne for
a time period of one year between April 2011- April 2012.}

\begin{abstract}
We study a mean value of the classical additive divisor problem, that is
\[
\sum_{f\sim F}\sum_{n\sim N}\left|\sum_{l\sim L}d(n+l)d(n+l+f)-\textrm{main term}
\right|^{2},
\]
with quantities $N\geq 1$, $1\leq F\ll N^{1-\varepsilon}$ and $1\leq L\leq N$. The main term we are interested in here is the one by Motohashi
\cite{Motohashi4}, but we also give an upper bound for the case where the main term is that of Atkinson \cite{Atkinson}.
Furthermore, we point out that the proof yields 
an analogous upper bound for a shifted convolution sum over Fourier coefficients of a fixed holomorphic cusp form in mean. 
\end{abstract}

\section{Introduction}

A widely studied occurrence of the classical
divisor function $d(n)=\sum_{d|n}1$ is in the additive divisor problem, in which one
investigates the asymptotic behavior of the sum
\begin{equation}\label{divisori}
D(x;m)=\sum_{n\leq x}d(n)d(n+m)
\end{equation}
as $x$ tends to infinity and $m$ is a given positive integer.
Ingham was the first to find an asymptotic formula for the above sum in 1927, and in 1931 Estermann brought the Kloosterman sums
into discussion, improving the earlier result. The importance of uniformity in the shift $m$ was noted ten years later by Atkinson \cite{Atkinson}, 
who further used his findings in studying the fourth power mean of the Riemann
zeta-function on the critical line. In 1982 Deshouillers and Iwaniec applied Kuznetsov's trace formulas into the problem, 
obtaining a notable improvement to the upper bound of the error term. Finally, especially crucial for
this paper, Motohashi \cite{Motohashi4} worked out an explicit spectral decomposition for a weighted convolution sum, separating a main term and obtaining an
upper bound for the error term uniformly in the shift $m$. We mention shortly papers by Heath-Brown, Iwaniec \cite{Iwaniec1} and Ivi\'{c}-Motohashi \cite{Ivic2}, \cite{Ivic},
who obtained mean value results for the error term, and \cite{Kuznetsov}, \cite{Meurman} and \cite{Tahtadjan} for other related results.
For a more thorough discussion of the history and the references of the above mentioned classical papers we refer to Motohashi's comprehensive paper \cite{Motohashi4}.

In recent years the analogy between the divisor function and the Fourier coefficients of holomorphic and non-holomorphic cusp
forms has appeared repeatedly in the literature. The behaviour  
of the analogous sums to \eqref{divisori},
\begin{equation}\label{holom}
A(x;m)=\sum_{n\leq x}a(n)\overline{a(n+m)}
\end{equation}
involving the Fourier coefficients of a holomorphic cusp form and
\begin{equation}\label{Maass}
T(x;m)=\sum_{n\leq x}t(n)t(n+m)
\end{equation}
involving the Hecke eigenvalues corresponding to Fourier coefficients
of a Maass form, has been studied intensively. In 1965 Selberg studied the meromorphity of a certain Dirichlet series
related to the Fourier coefficients of holomorphic cusp forms, while
Good in 1981 estimated the second moment of a modular $L$-function on the critical line using the spectral decomposition of shifted convolution sums 
coming from the Fourier coefficients of a holomorphic cusp form. Jutila \cite{Jutila1}, \cite{Jutila5} derived explicit asymptotic formulae for 
all three cases \eqref{divisori}, \eqref{holom} and \eqref{Maass} in a unified way, using the respective generating Dirichlet series. For a discussion of the substantial history including e.g. the works by 
Blomer, Duke, Friedlander, Harcos, Iwaniec, Jutila, Michel, Motohashi and Sarnak as well as a general spectral decomposition for the weighted shifted convolution sums we refer to the paper \cite{Harcos} by Blomer and Harcos. Along with these references we mention also \cite{BlomerHarcos}, \cite{Goldfeld}, \cite{Good3}, \cite{Jutila9} and \cite{Jutila6}.
Together the three sums above go under the name of the shifted convolution problem.

As a prologue for this paper, in Lemma 3 of his paper \cite{Jutila9} Jutila studies the sum \eqref{holom}
in mean, and his argument can easily be extended to prove
\begin{equation}\label{aasumma}
\sum_{0\leq f\leq F}\sum_{1\leq n\leq N}\bigg|\sum_{1\leq l\leq L}
a(n+l)\overline{a(n+f+l)}\bigg|^{2}\ll (N+F)^{k}N^{k}L
\end{equation}
for all $N,F\geq 1$ and $1\leq L\leq N$.
Instantaneously the proof gives also an upper bound 
for a sum over Hecke eigenvalues of Maass forms in mean, as stated in Lemma 6 of \cite{Jutila7} :  For $N,F\geq 1$ and $1\leq L\leq N$ 
\begin{equation}\label{Tsumma}
\sum_{0\leq f\leq F}\sum_{1\leq n\leq N}\bigg|\sum_{1\leq l\leq L}
t(n+l)t(n+f+l)\bigg|^ {2}\ll (N+F)^{1+\varepsilon}NL.
\end{equation}
The idea in the estimation of the triple sum lies
in its sensitiveness of the size of the innermost sum over $l$.
A similar result, sensitive for the size of the shift, is needed in the doctoral thesis of the author (\cite{Suvitie}, Lemma 3.4). 
An analogous result over more general settings also appears on p. 81 in the paper \cite{Blomer} by Blomer, Harcos and Michel, leading
to the same bounds \eqref{aasumma}, \eqref{Tsumma}.
In all the three papers \cite{Jutila9}, \cite{Jutila7}, \cite{Blomer} the motivation behind these studies has been gaining information about the 
upper bounds of the original shifted convolution sums \eqref{holom} and \eqref{Maass}, while the author needed her estimate in her thesis
in estimating a certain spectral sum over inner products involving holomorphic cusp forms and Maass forms.

However, to our knowledge there does not exist a similar bound for the case of the divisor function, and the extension of the proofs of the cuspidal
cases seem problematic due to the main term of the additive divisor problem. This paper is our 
attempt to bring new light for this case. Hence we state the following Theorem:
\begin{thm}[Main theorem]\label{thm1} 
Let $N\geq 1$, $1\leq L\leq N$ and $1\leq F\ll N^{1-\varepsilon}$. Then
\[
S=\sum_{f\sim F}\sum_{n\sim N}\left|\sum_{l\sim L}d(n+l)d(n+l+f)-\frac{6}{\pi^{2}}\int_{(L+n)/f}^{(2L+n)/f}m(x;f)dx
\right|^{2}
\]
\begin{equation}\label{paatermi}
\ll N^{2+\varepsilon}+N^{1+\varepsilon}LF.
\end{equation}
Here $m(x;f)$ is as in (1.12) in \cite{Motohashi4}.
\end{thm}
The notation $m\sim M$ stands for $M< m\leq2M$. 
For the classical main term $M(N,f)$ by Atkinson (\cite{Atkinson}, p. 185) we easily obtain the following upper bound:
\begin{cor}\label{cor2} 
Let $N\geq 1$, $1\leq L\leq N$ and $1\leq F\ll N^{1-\varepsilon}$. Then
\[
S=\sum_{f\sim F}\sum_{n\sim N}\left|\sum_{l\sim L}d(n+l)d(n+l+f)-M(2L+n,f)+M(L+n,f)
\right|^{2}
\]
\[
\ll N^{2+\varepsilon}+N^{1+\varepsilon}LF+N^{-1+\varepsilon}L^{2}F^{3}.
\]
\end{cor}

\begin{rem}
With small
values of the quantities $F$ and $L$, that is, if $L^{2}F=o(N)$, the trivial bound $N^{1+\varepsilon}L^{2}F$  yields
a better result than Theorem \ref{thm1} and Corollary \ref{cor2}. On the other hand, in their paper \cite{Ivic} Ivi\'{c} and Motohashi gain the upper bound
\[
\sum_{f=1}^{F}\left|\sum_{n=1}^{N}d(n)d(n+f)-\frac{6}{\pi^{2}}\int_{0}^{N/f}m(x;f)dx
\right|^{2}\ll N^{4/3+\varepsilon}F^{1/3}
\]
uniformly for $F\leq N^{1/2-\varepsilon}$. Here, according to Meurman (\cite{Meurman}, pp. 225-226)
it seems possible to extend the range to $F\leq N$ by replacing the upper bound by $N^{4/3+\varepsilon}F^{1/3}+NF$.
This result yields, up to factor $N^{\varepsilon}$, the same upper bound as in Theorem \ref{thm1} in case $F\gg N^{1/2}$ and $L=N$.
Other than this, to our knowledge, our result \eqref{paatermi} and the trivial bound together yield the genuinely best upper bound with the restriction $F\ll N^{1-\varepsilon}$.
\end{rem}

The main novelty in this paper is the methodological uniformity allowing 
an analogous treatment for the cuspidal cases:
Following the analogous proof for the holomorphic cusp forms we immediately gain
\begin{thm}\label{thm3} 
Let $N\geq 1$, $1\leq L\leq N$ and $1\leq F\ll N^{1-\varepsilon}$. Then
\[
\sum_{f\sim F}\sum_{n\sim N}\left|\sum_{l\sim L}a(n+l)\overline{a(n+l+f)}
\right|^{2}
\ll
N^{2k+\varepsilon}+N^{2k-1+\varepsilon}LF.
\]
\end{thm}

To compare our result to the earlier mentioned, we notice that up to a factor $N^{\varepsilon}$ and with the restriction 
$F\ll N^{1-\varepsilon}$, \eqref{aasumma} yields our result only if $L\asymp 1$, and
otherwise our result is genuinely better. 

In the case of non-holomorphic cusp forms we face the problem of the lack of a proper analogy for the spectral decomposition of the shifted convolution sum
in question. There exists an analogous result for our crucial Lemma \ref{Motohas}, namely Lemma 5 in \cite{Motohashi}, 
originally due to Jutila \cite{Jutila2}, \cite{Jutila3}. However, although the analogous part of this decomposition
results the anticipated analogous bound, there appears also an arithmetic correction term, which leads to additional problems. If we were dealing with
a genuinely oscillating weight function in our proof, then the contribution of the correction term would easily be
proven to stay under the expected bound. However, as this is not the case we shall lighten this paper and just assign the conjecture below,
leaving the proof for future study.
\begin{conj}\label{con} 
Let $N\geq 1$, $1\leq L\leq N$ and $1\leq F\ll N^{1-\varepsilon}$. Then
\[
\sum_{f\sim F}\sum_{n\sim N}\left|\sum_{l\sim L}t(n+l)t(n+l+f)
\right|^{2}
\ll
N^{2+\varepsilon}+N^{1+\varepsilon}LF.
\]
\end{conj}

\begin{rem}
Note also Theorem 1 in \cite{Harcos}, which gives a more general spectral decomposition of a shifted convolution sum including
Hecke eigenvalues of the underlying cusp forms of two arbitrary cuspidal automorphic representations of $PGL_{2}(\mathbb{Z})\backslash PGL_{2}(\mathbb{R})$.
\end{rem}

\section{The needed notation and auxiliary lemmas}

\subsection{The cusp forms}

Besides their analogous behavior to the divisor function, the holomorphic and non-holomorphic cusp 
forms play a role in the spectral decomposition of the shifted
convolution sum over the divisor function (Lemma \ref{Motohas}). Hence we start by introducing briefly some results concerning them.
For the proofs and for a general reference the reader is referred to Motohashi's monograph \cite{Motohashi2}.

We may confine ourselves to the cusp forms for the full modular group
$\Gamma=SL_{2}(\mathbb{Z})$ operating through M\"{o}bius
transformations on the upper half plane
$\mathbb{H}$. A
\textit{holomorphic cusp form}
$F(z):\mathbb{H}\rightarrow\mathbb{C}$ of weight $k\in\mathbb{Z}$
with respect to $\Gamma$ can be represented by its Fourier series
\[
F(z)=\sum^{\infty}_{n=1}a(n)e(nz), \;\; e(\alpha)=\exp(2\pi i
\alpha).
\]
We may assume that $k$ is even and $k\geq 12$, otherwise $F(z)$ is
trivial.
We let
\[
\{\psi_{j,k}\;|\; 1\leq j\leq\vartheta(k)\}
\]
be an orthonormal basis of the unitary space of holomorphic cusp
forms of weight $k$, and write
\[
\psi_{j,k}(z)=\sum_{n=1}^{\infty}\rho_{j,k}(n)n^{\frac{k-1}{2}}e(nz).
\]
We may suppose that the basis vectors are eigenfunctions of the
Hecke operators $T_{k}(n)$ for all positive integers $n$. Thus, in particular,
$T_{k}(n)\psi_{j,k}=t_{j,k}(n)\psi_{j,k}$ for certain real numbers
$t_{j,k}(n)$, which we call Hecke eigenvalues. Comparing Fourier
coefficients on both sides, one may verify that
$\rho_{j,k}(n)=\rho_{j,k}(1)t_{j,k}(n)$ for all $n\geq 1$, $1\leq
j\leq \vartheta(k)$.
We put $\alpha_{j,k}=16 \Gamma(2k)(4\pi)^{-2k-1}$ $\times|\rho_{j,k}(1)|^{2}$,
whence
\begin{equation}\label{thetaalpha}
\sum_{j=1}^{\vartheta(k)}\alpha_{j,k}\ll k.
\end{equation}
Furthermore, we let
\[
H_{j,k}(s)=\sum_{n=1}^{\infty}t_{j,k}(n)n^{-s}
\]
stand for the Hecke L-function attached to $\psi_{j,k}$. The
series converges absolutely for Re $s>1$ because of the bound
\begin{equation}\label{tee}
|t_{j,k}(n)|\leq d(n)\ll n^{\varepsilon}
\end{equation}
by Deligne \cite{Deligne}. Furthermore $H_{j,k}(s)$ can be continued 
to an entire function, and it satisfies a functional equation, which implies that for bounded $s$
\begin{equation}\label{hoojiikoo}
H_{j,k}(s)\ll k^{c}
\end{equation}
uniformly in $j$. Here $c$ is a suitable constant, which depends only on  Re $s$.

A \textit{non-holomorphic cusp form} $u(z)=u(x+iy):\mathbb{H}\rightarrow\mathbb{C}$ is a
non-constant real-analytic $\Gamma$-invariant function in the
upper half-plane, square-integrable with respect to the hyperbolic
measure $d\mu(z)=\frac{dx\,dy}{y^{2}}$ over a fundamental domain
of $\Gamma$. Furthermore $u(z)$ is an eigenfunction of the non-euclidean Laplacian
 $\Delta=-y^2(\frac{\partial^2}{\partial x^2}+\frac{\partial^2}{\partial
 y^2})$, and the corresponding eigenvalue can be written
as $1/4+\kappa^2$ with $\kappa>0$.
The Fourier series expansion for $u(z)$ is then of the form
\[
u(z)=y^{1/2}\sum_{n\neq 0}\rho(n)K_{i\kappa}(2\pi |n|y)e(nx)
\]
with $K_{\nu}$ a Bessel function of imaginary argument. We may
suppose that our cusp forms are eigenfunctions of the Hecke
operators $T(n)$ for all positive integers $n$ and that $u(x+iy)$
is even or odd as a function of $x$. Thus $T(n)u=t(n)u$ for
certain real numbers $t(n)$, which are again called Hecke eigenvalues,
and $u(-\overline{z})=\pm u(z)$. Comparing Fourier coefficients on
both sides, one may verify that $\rho(n)=\rho(1)t(n)$ and
$\rho(-n)=\pm\rho(n)$ for all $n\geq 1$. 
The \textit{Maass (wave) forms} $u_{j}$
constitute an orthonormal set of non-holomorphic cusp forms
arranged so that the corresponding parameters $\kappa_{j}$
determined by the eigenvalues $1/4+\kappa_{j}^{2}$ lie in an
increasing order. We write $\rho_{j}(n)$ and $t_{j}(n)$ for the corresponding Fourier coefficients and Hecke
eigenvalues.
Now let
\[
\alpha_{j}=|\rho_{j}(1)|^{2}/\cosh(\pi\kappa_{j}),
\]
and let
\[
H_{j}(s)=\sum_{n=1}^{\infty}t_{j}(n)n^{-s}
\]
stand for the Hecke $L$-function attached to $u_{j}$. As the Ramanujan-Petersson conjecture
$t_{j}(n)\ll n^{\varepsilon}$ holds in a mean value sense, the series converges
absolutely for Re $s>1$. Again, $H_{j}(s)$ can be continued to an entire function.
We have the bound
\begin{equation}\label{hoonelja}
\sum_{\kappa_{j}\leq K}\alpha_{j}H_{j}^{4}\left(\frac{1}{2}\right)\ll K^{2}\log^{15}K.
\end{equation}
(See \cite{Motohashi2}, Theorem 3.4.)

Lastly, in this paper, the following notation will be adopted:  
Vinogradov's relation $f(z)\ll g(z)$ is another notation for
$f(z)=\mathcal{O}(g(z))$. We let $\varepsilon$ stand generally for a small positive number,
not necessarily the same at each occurrence. 

\subsection{The needed lemmas}

In this section we shall gather some auxiliary results which will be used
during the course of the proofs of the theorems.

For estimating the
Gamma function we shall frequently use the following 
Stirling's formula without further reference:
In any fixed strip $b\leq\sigma\leq c$ we have
\[
|\Gamma(\sigma+i
t)|=\sqrt{2\pi}|t|^{\sigma-1/2}e^{-|t|\pi/2}(1+\mathcal{O}(|t|^{-1})),
\]
for $|t|\rightarrow \infty$.
For a proof, see Olver \cite{Olver} p. 294.
Further we have 
\begin{equation}\label{kertoma}
\Gamma(x)= \sqrt{2\pi}x^{x-1/2}e^{-x}(1+\mathcal{O}(x^{-1})),
\end{equation}
for all positive real $x\rightarrow \infty$, see \cite{Lebedev}, Eq. (1.4.25).

The proof of the following classical estimate for the fourth moment of Riemann's zeta-function on the critical line
can be found e.g. from \cite{Titchmarsh}, p. 147:
\begin{equation}\label{zetalog}
\int_{0}^{T}\left|\zeta\left(\frac{1}{2}+it\right)\right|^{4}\,dt\ll T \log^{4} T
\end{equation}
for $T\geq 1$. Along with this estimate we have
\begin{equation}\label{zeta12}
\int_{0}^{T}\left|\zeta\left(\frac{1}{2}+it\right)\right|^{12}\,dt\ll T^{2} \log^{17} T
\end{equation}
for $T\geq 1$ by Heath-Brown \cite{Heath-Brown}.
Moreover it is known that
\[
\frac{1}{\zeta(1+it)}\ll \log |t|,
\]
as $|t|\geq 1$. For a proof, see e.g. 
\cite{Titchmarsh}, p. 132.

We
have an important tool arising from spectral theory:
\begin{lem}[The spectral large sieve]
For $K\geq 1$, $1\leq \Delta\leq K$, $M\geq 1$ and any complex
numbers $a_{m}$ we have
\[
\sum_{K\leq\kappa_{j}\leq K+\Delta}\alpha_{j}\left|\sum_{m\leq
M}a_{m}t_{j}(m)\right|^{2}\ll
(K\Delta+M)(KM)^{\varepsilon}\sum_{m\leq M}|a_{m}|^{2}.
\]
\end{lem}
For a proof, see Theorem 1.1 in \cite{Jutila4} or Theorem 3.3 in
\cite{Motohashi2}.
The continuous analogy is the following:
\begin{lem}
For $K$ real, $\Delta\geq 0$, $M\geq 1$ and any complex
numbers $a_{m}$ we have
\[
\int_{K}^{K+\Delta}\left|\sum_{m\leq
M}a_{m}\sigma_{2ir}(m)m^{-ir}\right|^{2}\, dr\ll
(\Delta^{2}+M)M^{\varepsilon}\sum_{m\leq M}|a_{m}|^{2}
\]
uniformly in $K$, as $\sigma_{2ir}(f)=\sum_{d|f}d^{2ir}$.
\end{lem}
Proof can be found in \cite{Suvitie}, Lemma 1.12.

A crucial role is played by a spectral decomposition of the
shifted convolution sum over the divisor function:
\begin{lem}\label{Motohas}
Let $f$ be a positive integer and $W$ a smooth function of compact
support on $(0,\infty)$. Then
\[
\sum_{l=1}^{\infty}d(l)d(l+f)W\left(\frac{l}{f}\right)=\frac{6}{\pi^{2}}\int_{0}^{\infty}
m(x;f)W(x)\, dx
+\frac{f^{1/2}}{\pi}\int_{-\infty}^{\infty}f^{-ir}\sigma_{2ir}(f)
\]
\[
\times \frac{|\zeta(\frac{1}{2}+ir)|^{4}}{|\zeta(1+2ir)|^{2}}\Theta(r;W)\,dr
+f^{1/2}\sum_{j=1}^{\infty}\alpha_{j}t_{j}(f)H_{j}^{2}\left(\frac{1}{2}\right)\Theta(\kappa_{j};W)
\]
\begin{equation}\label{spekthajo}
+\frac{1}{4}f^{1/2}\sum_{k=6}^{\infty}\sum_{j=1}^{\vartheta(k)}(-1)^{k}\alpha_{j,k}t_{j,k}(f)H_{j,k}^{2}\left(\frac{1}{2}\right)\Xi\left(k-\frac{1}{2};W\right)
\end{equation}
Here $m(x;f)$ is defined by (1.12) in \cite{Motohashi4}, $\sigma_{2ir}(f)$ is as in the previous lemma
and
\[
\Theta(r;W) =\frac{1}{2}\textrm{Re}\left\{\left(1+\frac{i}{\sinh(\pi r)}\right)\Xi(ir;W)\right\}
\]
\begin{equation}\label{theeta}
=\frac{1}{4}\left(1+\frac{i}{\sinh(\pi r)}\right)\Xi(ir;W)+(r\mapsto -r),
\end{equation}
with
\[
\Xi(\xi;W)=\frac{\Gamma(\xi+\frac{1}{2})^{2}}{\Gamma(2\xi+1)}
\int_{0}^{\infty}x^{-1/2-\xi}
F\left(\xi+\frac{1}{2},\xi+\frac{1}{2};2\xi+1;-\frac{1}{x}\right)
\]
\begin{equation}\label{ksiii}
\times W(x)
\,dx.
\end{equation}
Here $F(. \, , .\, ; . \, ; .)$ is the hypergeometric function and $(r\mapsto -r)$ stands for an expression similar to the
preceding one, but with $r$ replaced by $-r$.
\end{lem}
For the proof, see \cite{Motohashi4}, Theorem 3.

Finally we introduce the basic inequality in the proof of the
classical large sieve:
\begin{lem}[Sobolev]\label{klassinen}
Let $a\leq u\leq a+\Delta$ and let the function $f$ be
continuously differentiable on this interval. Then
\[
|f(u)|^{2}\leq \Delta^{-1}\int_{a}^{a+\Delta}|f(x)|^{2}\,dx
+2\left(\int_{a}^{a+\Delta}|f(x)|^{2}\,dx\right)^{1/2}
\]
\[
\times \left(\int_{a}^{a+\Delta}|f'(x)|^{2}\,dx\right)^{1/2} \ll
\Delta^{-1}\int_{a}^{a+\Delta}\left(|f(x)|^{2}+\Delta^{2}|f'(x)|^{2}
\right)\,dx
\]
uniformly.
\end{lem}
For a proof, see Montgomery \cite{Montgomery}, Lemma 1.1 applied
to $f^{2}$.

\section{Proofs of the results}

\subsection{Theorem \ref{thm1}}

We start by adding a set of $\asymp\log(1/U)$ real-valued smooth weight functions
$g_{\delta}\left(\frac{l}{L}\right)$ to the $l$-sum so that
their sum produces an approximation of the characteristic function of the interval $[1,2]$ with an error of size $\ll U$
and their supports 
widen step by step by factors 2 when we move away from the end points 1 and 
2. Thus if the supports of the first and last weight functions are of length $U=L^{-1+\varepsilon}$, then the next 
ones are of length $\asymp 2U$ and so on, and the slopes of the weight functions cancel out each other. 
Hence we let $g_{\delta}(x)$ stand for a compactly supported function on some interval of length $\asymp\delta$ contained in $[1,2]$. Moreover,  $g(x)=1$ on
an interval of length $\asymp\delta$ and
$g^{(\nu)}(x)\ll_{\nu}\delta^{-\nu}$ for each $\nu\geq 0$ and $x\in \mathbb{R}$. Here $L^{-1+\varepsilon}\leq\delta\ll 1$. Now
\[
S\ll L^{\varepsilon}\sum_{f\sim F}\sum_{n\sim N}\bigg |\sum_{l=1}^{\infty}d(n+l)d(n+l+f)W\left(\frac{n+l}{f}\right)
\]
\[
-\frac{6}{\pi^{2}}\int_{0}^{\infty}m(x;f)
 W(x)dx\bigg|^{2}+N^{1+\varepsilon}FU^{2}L^{2}
\]
\[
\ll L^{\varepsilon}\sum_{f\sim F}\sum_{n\sim N}\left| e_{1}(n,f;L)+ e_{2}(n,f;L)+ e_{3}(n,f;L)\right|^{2}
+N^{1+\varepsilon}F,
\]
with $W(x)=g_{\delta}\left(\frac{fx-n}{L}\right)$. Here $e_{\nu}(n,f;L)$ ($\nu=1,2,3$) stands for the $(\nu+1)$th term on the right hand side of \eqref{spekthajo}.

The contribution of the term  $e_{2}(n,f;L)$ turns out to be the most difficult to bound. Let us denote
\[
S_{2}=\sum_{f\sim F}\sum_{n\sim N}\left| f^{1/2}\sum_{j=1}^{\infty}\alpha_{j}t_{j}(f)H_{j}^{2}\left(\frac{1}{2}\right)\Theta(\kappa_{j};W)\right|^{2}
\]
with $\Theta(r;W)$ as in \eqref{theeta}. As the support of $W$ tends to
infinity, we have in the integral in \eqref{ksiii}
\[
F\left(\frac{1}{2}+ir,\frac{1}{2}+ir;1+2ir;-\frac{1}{x}\right)
=
\left(\frac{1+\sqrt{1+x^{-1}}}{2}\right)^{-1-2ir}
\]
\begin{equation}\label{hyperg}
\times F\left(\frac{1}{2}+ir,\frac{1}{2};1+ir;\left(\frac{1-\sqrt{1+x^{-1}}}{1+\sqrt{1+x^{-1}}}\right)^{2}\right)
\end{equation}
by the quadratic transformation formula for the hypergeometric
function from Lebedev \cite{Lebedev}, Eq. (9.6.12). We write the RHS of \eqref{hyperg} by the definition of the hypergeometric series. 
Now clearly the $k$th term in the series is non-oscillating and of order $x^{-2k}$. Therefore it suffices to study
only the leading term $1$, others behaving similarly until some sufficient constant, and the contribution of the rest being negligibly small.
(See also Convention 2 in \cite{JutilaMotohashi}.)

By integrating repeatedly by parts we notice that $\kappa_{j}$ can be
truncated to be $\ll N^{1+\varepsilon}\delta^{-1}L^{-1}$ with negligibly small error. We further decompose the range
$[1,N^{1+\varepsilon}\delta^{-1}L^{-1}]$ into dyadic intervals $\kappa_{j}\sim K$ with $1\leq K\ll N^{1+\varepsilon}\delta^{-1}L^{-1}$, the number
of these intervals being $\ll N^{\varepsilon}$.
Next we prepare ourselves to use the duality principle (see for example pp. 169-170 in \cite{Kowalski}). 
Let $b_{f,n}$ stand for any finite sequence of complex numbers such that $\sum_{f\sim F}\sum_{n\sim N}|b_{f,n}|^{2}=1$.
Then
\[
S_{2}\ll N^{\varepsilon}L^{2}F^{-1}\sum_{\kappa_{j}\sim K}\alpha_{j}\left| \sum_{f\sim F}\sum_{n\sim N}b_{f,n}t_{j}(f)I(\kappa_{j};W)\right|^{2}
\sum_{\kappa_{j}\sim K}\alpha_{j}
H_{j}^{4}\left(\frac{1}{2}\right)\kappa_{j}^{-1}
\]
with
\[
I(\kappa_{j};W)=\int_{1}^{2}\left(\frac{zL+n}{f}\right)^{-1/2-i\kappa_{j}}
\left(\frac{1+\sqrt{1+\frac{f}{zL+n}}}{2}\right)^{-1-2i\kappa_{j}}g_{\delta}(z)\, dz.
\]
Note that  $\kappa_{j}$ replaced by $-\kappa_{j}$ can be taken care of with complex conjugation.
By \eqref{hoonelja} this is further
\begin{equation}\label{valiin}
\ll  N^{\varepsilon}KL^{2}F^{-1}\sum_{\kappa_{j}\sim K}\alpha_{j}\left| \sum_{f\sim F}\sum_{n\sim N}b_{f,n}t_{j}(f)I(\kappa_{j};W)\right|^{2}.
\end{equation}

We want to next apply the spectral large sieve, so we shall use Sobolev's lemma \ref{klassinen} as done in Jutila's paper \cite{Jutila1} on p. 454 to get rid of 
the parameter
$\kappa_{j}$ in the integral $I(\kappa_{j};W)$. Therefore  
the range $(K,2K]$ for $\kappa_{j}$ is
split up into segments of length $\Delta$ in such a way that the
integral $I(y;W)$ remains essentially stationary as $y$ runs
over a segment. That is
\[
\Delta \log\left(\frac{(zL+n)\left(1+\sqrt{1+\frac{f}{zL+n}}\right)^{2}}{4f}\right)\ll \log N.
\]
In this way, the second term in the upper bound in Lemma
\ref{klassinen} will be comparable to the first. An appropriate choice would be
$\Delta=1$.

Now we divide the $\kappa_{j}$-sum in \eqref{valiin} into
subsums of length $\Delta$, and apply Lemma \ref{klassinen} to each subsum
arriving at the bound
\[
N^{\varepsilon}KL^{2}F^{-1}\sum_{l=0}^{
K-1}\int_{K+l}^{K+l+1}
\sum_{K+l<\kappa_{j}\leq K+l+1}\alpha_{j}\bigg(\bigg| \sum_{f\sim F}\sum_{n\sim N}b_{f,n}t_{j}(f)I(y;W)\bigg|^{2}
\]
\[
+\bigg| \sum_{f\sim F}\sum_{n\sim N}b_{f,n}t_{j}(f)\frac{\partial}{\partial y}I(y;W)\bigg|^{2}\bigg)\, dy.
\]

We next apply the spectral large sieve to the subsum over
$\kappa_{j}$, and finally add the results
corresponding to all subsums. This leads us to the bound 
\[
N^{\varepsilon}KL^{2}F^{-1}(K+F)
\sum_{f\sim F}\int_{K}^{2K}\bigg(\bigg|\sum_{n\sim N}b_{f,n}I(y;W)\bigg|^{2}
\]
\[
+\bigg| \sum_{n\sim N}b_{f,n}\frac{\partial}{\partial y}I(y;W)\bigg|^{2}\bigg)\, dy.
\]
Now we use Cauchy's inequality for the integrals $I(y;W)$ and $\frac{\partial}{\partial y}I(y;W)$, and attach
a suitable real-valued smooth weight function
$u(y)$ to the $y$-integral. We set $u$ to be compactly supported on the interval $[K/2,3K]$, $u(y)=1$ on
$[K,2K]$ and $u^{(\nu)}(y)\ll_{\nu}K^{-\nu}$ for each $\nu\geq 0$ and $y\in \mathbb{R}$. By opening the squares
we end up with a double sum $\sum_{n_{1}\sim N}\sum_{n_{2}\sim N}$. By repeated partial integration over the $y$-integral we now
see that we may truncate $n_{1}-n_{2}\ll N^{1+\varepsilon}K^{-1}$. Hence we finally have an upper bound
\[
N^{-1+\varepsilon}K^{2}L^{2}(K+F)\delta^{2}
\sum_{f\sim F}\sum_{n\sim N}|b_{f,n}|^{2}\sum_{|n-n'|\ll N^{1+\varepsilon}K^{-1}}1,
\]
and we conclude with the desired upper bound.

The contribution of $e_{1}(n,f;L)$ to $S$ can be estimated by exactly the same argumentation as above. We need the estimation
\begin{equation}\label{zeta8}
\int_{0}^{T}\left|\zeta\left(\frac{1}{2}+it\right)\right|^{8}\,dt\ll T^{3/2+\varepsilon} 
\end{equation}
following readily of \eqref{zetalog} and \eqref{zeta12} by the Cauchy's inequality. Hence
\[
S_{1}\ll N^{1+\varepsilon}LF.
\]
We could use also a more straightforward argumentation, but these steps make it easier to follow the proof of Theorem \ref{thm3}.

Using \eqref{thetaalpha}, \eqref{tee}, \eqref{hoojiikoo} and \eqref{kertoma} the trivial estimates suffice for the estimation of $e_{3}(n,f;L)$.

\subsection{Corollary \ref{cor2}}

Using straightforward calculation it is easily seen that the difference of our main term and that of Atkinson's is
\[
\frac{6}{\pi^{2}}\int_{(L+n)/f}^{(2L+n)/f}m(x;f)dx-
(M(2L+n,f)-M(L+n,f))\ll FLN^{-1+\varepsilon},
\]
whence the result follows.

\subsection{Theorem \ref{thm3}}

Analogous to the proof of Theorem \ref{thm1}. For the size of the Fourier coefficients see \eqref{tee}. Instead of Lemma \ref{Motohas} we use Lemma 4 in \cite{Motohashi},
originally due to Motohashi \cite{Motohashi5}. The counterpart for \eqref{hoonelja} is Lemma 4 in \cite{Jutila1}, and instead of \eqref{zeta8} we need
Lemma 3 in \cite{Jutila1}.

\subsection*{Acknowledgment}
The author wishes to thank Professor Matti Jutila for the idea for this
subject as well as for pointing out the advantage of choosing the weight function in the proof of Theorem \ref{thm1} as it was done, leading to
considerable improvement of the results. His help in finding a proof of the duality principle easily extendable to the continuous case was
of great value. The kind support and advices of Gergely Harcos, Tom Meurman, Philippe Michel and Yoichi Motohashi are gratefully acknowledged.

\vspace{1cm}

{\sc Eeva Suvitie}

{\sc Department of Mathematics}

{\sc FI-20014 University of Turku} 

{\sc Finland}

{\sc e-mail: eeva.suvitie@utu.fi}

\end{document}